\documentclass{amsart}

\numberwithin{equation}{section}
\newcommand{\defi}{\stackrel{\textup{\tiny def}}{=}}

\DeclareMathOperator{\spa}{span} \DeclareMathOperator{\ran}{ran}

\newtheorem{theorem}{Theorem}[section]
\newtheorem{lemma}[theorem]{Lemma}
\newtheorem{proposition}[theorem]{Proposition}

\theoremstyle{definition}
\newtheorem{definition}[theorem]{Definition}

\theoremstyle{remark}
\newtheorem{remark}[theorem]{Remark}
\begin{document}

\title[Carath\'eodory functions]{Carath\'eodory
functions in the Banach space setting}

\author[D. Alpay]{Daniel Alpay}
\address{Department of Mathematics\\
Ben--Gurion University of the Negev\\
Beer-Sheva 84105\\ Israel} \email{dany@math.bgu.ac.il}
\thanks{D. Alpay thanks the Earl Katz family for endowing the chair which
supports his research}

\author[O. Timoshenko]{Olga Timoshenko}
\address{Department of Mathematics\\
Ben--Gurion University of the Negev\\
Beer-Sheva 84105\\ Israel} \email{olgat@math.bgu.ac.il}

\author[D. Volok]{Dan Volok}
\address{Department of Mathematics\\
Kansas State University\\
Manhattan, KS 66506-2602\\USA}\email{danvolok@hotmail.com}

\date{}

\begin{abstract}
We prove representation theorems for Carath\'eodory functions in the setting of Banach spaces.
\end{abstract}
\maketitle

\section{Introduction}
L. de Branges and J. Rovnyak introduced in \cite{dbr1}, \cite{dbr2}
various families of reproducing kernel Hilbert spaces of functions which take
values in a Hilbert space and are analytic in the open unit disk
or in the open upper half-plane. These spaces play an important role
 in operator theory, interpolation theory,
inverse scattering, the theory of wide sense stationary stochastic
processes and related topics; see for instance \cite{Dmk},
\cite{ad1}, \cite{ad2}. In the case of the open unit disk
$\mathbb{D}$, of particular importance are the following two
kinds of reproducing kernels:
\begin{align}
\label{avron} k_\phi(z,w)&=\frac{\phi(z)+\phi(w)^*}{2(1-zw^*)},\\
\nonumber
k_s(z,w)&=\frac{I-s(z)s(w)^*}{1-zw^*}.
\end{align}
In these expressions, $s(z)$ and $\phi(z)$ are operator-valued functions analytic
in $\mathbb{D}$, $*$ denotes
the Hilbert space adjoint, and $I$ denotes the identity operator.
The functions for which the kernels $k_\phi(z,w)$ and $k_s(z,w)$ are positive
are called respectively Carath\'eodory and Schur functions.
We remark that one can use the Cayley
transform
$$s(z)=(I-\phi(z))(I+\phi(z))^{-1}$$
to reduce the study of the kernels $k_\phi(z,w)$ to the study of
the kernels $k_s(z,w)$. For these latter it is
well known that the positivity of the kernel $k_s(z,w)$ implies analyticity of $s(z)$.\\

Every Carath\'eodory function
admits two equivalent representations. The first, called the
Riesz -- Herglotz representation, reads as follows:
\begin{equation}
\label{anaelle}
\phi(z)=ia+\int_0^{2\pi}\frac{e^{it}+z}{e^{it}-z}d\mu(t)
\end{equation}
where $a$ is a real number and where $\mu(t)$ is an increasing
function such that $\mu(2\pi)<\infty$. The integral is a
Stieltjes integral and the proof relies on Helly's theorem; see
\cite[pp. 19--27]{MR48:904}.\\

The second representation reads:
\begin{equation}
\label{anaelle2}
\phi(z)=ia+\Gamma(U+zI)(U-zI)^{-1}\Gamma^*
\end{equation}
where $a\in{\mathbb R}$ and where $U$ is a unitary operator in an
auxiliary Hilbert space ${\mathcal H}$ and $\Gamma$ is a bounded
operator from ${\mathcal H}$ into ${\mathbb C}$.\\

The  expression \eqref{anaelle2} still makes sense in a more
general setting when  the kernel $k_\phi(z,w)$ has a finite
number of negative squares. The space ${\mathcal H}$ is then a
Pontryagin space. This is the setting in which Kre\u\i n and
Langer proved this result; see \cite[Satz 2.2 p. 361]{kl1}. They
allowed the values of the function $\phi(z)$ to be operators
between Pontryagin spaces and required weak continuity at the
origin. Without this hypothesis one can find functions for which
the kernel $k_\phi(z,w)$ has a finite number of negative squares
but which are not meromorphic in ${\mathbb D}$ and in particular
cannot admit representations of the form \eqref{anaelle2}; for
instance the function
$$\phi(z)=\begin{cases} 0\,\,{\rm if}\,\, z\not=0\\
                               1\,\,{\rm if}\,\, z=0
                               \end{cases}
                               $$
defines a kernel $k_\phi(z,w)$ which has one negative square; see
\cite[p. 82]{adrs} for an analogue for $k_s(z,w)$
kernels.\\

Operator--valued Carath\'eodory functions were extensively studied in the
Hilbert space case; see e.g. \cite{MR1821917}, \cite{MR1759548},
\cite{MR899274}, \cite{MR1036844}. We would also like to mention the
non-stationary setting, where an analogue of the representation
\eqref{anaelle2}  was obtained for  upper-triangular operators; see
\cite{MR2003f:47021} and \cite{MR1704663}.\\

The notion  of reproducing kernel space (with positive or
indefinite metric) can also be introduced for functions which
take values in  Banach spaces and even topological vector spaces.
The positive case was studied already by Pedrick for functions
with values in certain topological vector spaces in an
unpublished report \cite{pedrick} and studied further by P. Masani in his 1978 paper
\cite{masani78}. Motivations originate from the theory of partial
differential equations (see e.g. \cite{MR1973084}) and the theory
of stochastic processes (see e.g. \cite[\S 4]{MR647140} and \cite{MR626346}).\\

The present paper is devoted to the study of Carath\'eodory functions
whose values are bounded operators between appropriate Banach
spaces. It seems that there are no natural analogs of
Schur functions or of the Cayley transform in this setting.\\

 Let ${\mathcal B}$ be a Banach space. We denote by ${\mathcal B}^*$ the
 space of anti-linear bounded
functionals (that is, its conjugate dual space). The duality between ${\mathcal B}$
and ${\mathcal B}^*$ is
denoted by
\[
\langle b_*, b\rangle_{\mathcal B}\defi b_*(b),\quad \text{where
} b\in {\mathcal B}\text{ and } b_*\in{\mathcal B}^*.\] An
${\mathbf L}({\mathcal B},{\mathcal B}^*)$-valued function $\phi(z)$ defined in
some open neighborhood $\Omega$ of the origin
and weakly continuous at the origin will be called a Carath\'eodory
function if the ${\mathbf L}({\mathcal B},{\mathcal B}^*)$-valued
kernel
\begin{equation}
\label{toto}
k_\phi(z,w)=\frac{\phi(z)+\phi(w)^*\big|_{\mathcal B}}{2(1-zw^*)}
\end{equation}
is positive in $\Omega$. The notion of positivity for bounded
operators and kernels from ${\mathcal B}$
into ${\mathcal
B}^*$ is reviewed in the next section.  We shall prove (see Theorem \ref{Ivry-sur-Seine})
that every Carath\'eodory function admits a representation of the form
\eqref{anaelle2} and, in particular, admits an analytic extension to
 $\mathbb{D}$; see e.g. \cite[pp. 189--190]{MR58:12429a} for
information on vector-valued analytic functions. We note that the proof of
this theorem can be adapted to the case when the kernel
$k_\phi(z,w)$ has a finite number of negative squares;
see Remark \ref{indef}.\\

Furthermore, if $\mathcal{B}$ is a separable Banach space
then ${\mathbf L}({\mathcal B},{\mathcal
B}^*)$-valued  Carath\'eodory functions can be characterized as
 functions analytic in $\mathbb{D}$ and such that
\[
\phi(z)+\phi(z)^*\bigm|_{\mathcal B}\ge 0,\quad z\in{\mathbb D}.
\]
Moreover, in this case we have an analogue of the Riesz -- Herglotz
representation \eqref{anaelle}; see Theorem
\ref{Antony}.\\

We conclude with the outline of the paper; the next three sections
are of preliminary nature, and deal with positive operators, Stieltjes
integrals and Helly's theorem respectively. Representation theorems for
Carath\'eodory functions are proved in Section 5. Two cases are to be
distinguished, as whether ${\mathcal B}$ is separable or not.
The case of  ${\mathbf L}({\mathcal B}^*,{\mathcal B})$-valued
Carath\'eodory functions will be treated in the last section of this paper.
This case is of special importance. Indeed, if $\phi(z)$ is a  ${\mathbf L}
({\mathcal B}^*, {\mathcal B})$-valued Carath\'eodory function which takes
invertible values, its inverse is a
${\mathbf L}({\mathcal B},{\mathcal B}^*)$-valued  Carath\'eodory
function. In the Hilbert space case, this fact has important connections
with operator models for pairs of unitary
operators (see \cite{dbs} and, for the
analogue for self--adjoint operators,  \cite{dbr1}, \cite{MR1960423}).
We will explore the Banach space generalizations of these results
in a future publication.

\section{Positive operators and positive kernels}
 In this section we review for the
convenience of the reader various facts on bounded positive operators from
${\mathcal B}$ into ${\mathcal B}^*$.
First some notations and a definition.

\begin{definition} Let  ${\mathcal B}$ be a  complex Banach space and let
$A\in{\mathcal L}({\mathcal B},{\mathcal B}^*).$ The operator $A$ is said to
be {\em positive} if
\[
\langle Ab,b\rangle_{\mathcal B}\geq 0,\quad \forall b\in {\mathcal B}.\]
\end{definition}
Note that a positive operator is in particular self-adjoint in the sense
that $A=A^*\bigm|_{\mathcal B}$, that
is,
\begin{equation}
\label{gabriel} \langle Ab,c\rangle_{{\mathcal B}}=
\overline{\langle Ac,b\rangle_{{\mathcal B}}}.
\end{equation}
Indeed, \eqref{gabriel} holds for $b=c$ in view of the positivity. It then
holds for all choices of $b,c\in{\mathcal B}$ by polarization:
\[
\begin{split}
\langle Ab,c\rangle_{{\mathcal B}}&=\frac{1}{4} \sum_{k=0}^3i^k\langle A(b+i^kc),(b+i^kc)
\rangle_{{\mathcal B}}\\
&=\frac{1}{4}
\sum_{k=0}^3i^k\langle A(c+i^{-k}b),(c+i^{-k}b)
\rangle_{{\mathcal B}}\\
&=\frac{1}{4}\overline{
\sum_{k=0}^3i^{-k}\langle A(c+i^{-k}b),(c+i^{-k}b)
\rangle_{{\mathcal B}}}\\
&=\overline{\langle Ac,b\rangle_{{\mathcal B}}}.
\end{split}
\]
Now, let $\tau$ be the natural injection
from ${\mathcal B}$ into
${\mathcal B}^{**}$:
\begin{equation}
\label{bastille}
\langle \tau(b), b_*\rangle_{{\mathcal B}^*} =
\overline{\langle b_*, b\rangle_{\mathcal B}}.
\end{equation}
We have for $b,c\in{\mathcal B}$:
\[
\begin{split}
\langle A^*\tau c,b\rangle_{\mathcal B}&=\langle \tau c, Ab
\rangle_{{\mathcal B}^*}\\
&=\overline{\langle Ab,c\rangle_{\mathcal B}}\\
&=\langle Ac,b\rangle_{\mathcal B}
\end{split}
\]
in view of \eqref{gabriel}, and hence $A=A^*\big|_{\mathcal B}$.\\

The following factorization result is well known and originates
with the works of Pedrick \cite{pedrick} (in the case of
topological vector spaces with appropriate properties) and
Vakhania \cite[\S 4.3.2 p.101]{MR626346} (for positive elements
in ${\mathbf L}({\mathcal B}^*,{\mathcal B})$); see the
discussion in \cite[p. 416]{masani78}. We refer also to
\cite{MR647140} for the case of barreled spaces and to
\cite{MR1973084} for the case of unbounded operators.

\begin{theorem}
The operator $A\in{\mathcal L}({\mathcal B},{\mathcal B}^*)$ is
positive if and only if there exist a Hilbert space ${\mathcal H}$
and a bounded operator $T\in{\mathbf L} ({\mathcal B},{\mathcal
H})$ such that $A=T^*T$. Moreover,
\begin{equation}
\label{samantha}
\langle Ab,c\rangle_{{\mathcal B}}=\langle Tb,
Tc \rangle_{\mathcal H},\quad b,c\in {\mathcal B}
\end{equation}
and
\begin{equation}
\label{bagdad} \sup_{\|b\|=1}\langle Ab,b\rangle_{{\mathcal B}} =\|A\|=\|T\|^2.
\end{equation}
Finally we have
\begin{equation}
\label{sof-sof}
|\langle Ab,c\rangle_{\mathcal B}|\le\langle Ab,b\rangle^{1/2}
\langle Ac,c\rangle^{1/2}.
\end{equation}
\label{rosh-hashana}
\end{theorem}
\begin{proof} For $Ab$ and $Ac$ in the range of $A$ the expression
\begin{equation}
\label{mind your step} \langle Ab,Ac\rangle_A= \langle Ab,c\rangle_{{\mathcal B}}=
\overline{\langle
Ac,b\rangle_{{\mathcal B}}}.
\end{equation}
is well defined in the sense that $Ab=0$ (resp. $Ac=0$) implies that
\eqref{mind your step} is equal to $0$.
Thus formula \eqref{mind your step} defines a 
sesquilinear form on $\ran A$. It is
positive since the operator $A$ is positive. Moreover, it is non-degenerate because
  if $\langle Ab,b\rangle_{\mathcal{B}}=0$ then $Ab=0.$ \\

Indeed, if $c$ is such that $\langle Ac,c\rangle_{\mathcal{B}}=0$ then 
\eqref{gabriel} implies that 
the real and the imaginary parts of $\langle Ab,c\rangle_{\mathcal{B}}$ are equal, respectively, to
\[\dfrac{1}{2}\langle A(b+c),b+c\rangle_{\mathcal{B}}\quad\text{and}\quad\dfrac{1}{2}\langle A(b+ic),b+ic\rangle_{\mathcal{B}}\] and, therefore, are non-negative. Then the same can be said about
$\langle Ab,-c\rangle_{\mathcal{B}},$ hence
\begin{equation}\label{referee}
\langle Ab,c\rangle_{\mathcal{B}}=0.
\end{equation}
Furthermore, if $c$ is such that $\langle Ac,c\rangle_{\mathcal{B}}>0$ then 
we have \[0\leq \left\langle A\left(b-\dfrac{\langle Ab,c\rangle_{\mathcal{B}}}{\langle Ac,c\rangle_{\mathcal{B}}}c\right),b-\dfrac{\langle Ab,c\rangle_{\mathcal{B}}}{\langle Ac,c\rangle_{\mathcal{B}}}c\right\rangle_{\mathcal{B}}=-\dfrac{|\langle Ab,c\rangle_{\mathcal{B}}|^2}{\langle Ac,c\rangle_{\mathcal{B}}},\]
hence \eqref{referee} holds in this case, as well.\\

  Thus, $(\ran A,\langle\cdot,\cdot\rangle_A)$ is a
pre-Hilbert
space. We will denote by ${\mathcal H}_A$ its completion and define
$$T\,:\,{\mathcal B}\longrightarrow {\mathcal H}_A,\quad Tb\defi Ab.$$
We have for $b,c\in{\mathcal B}$
\[
\langle T^*(Ac),b\rangle_{{\mathcal B}}= \langle Ac,Tb\rangle_{{\mathcal H}_A} =
\langle Ac,b\rangle_{{\mathcal
B}}.
\]
Hence, $T^*$ extends continuously to the injection map from ${\mathcal H}_A$
 into ${\mathcal B}^*$. We note that
\[
\langle Tb,Tc\rangle_{{\mathcal H}_A}=\langle T^*Tb,c\rangle_{{\mathcal B}}
=\langle Ab,c\rangle_{{\mathcal B}}\]
(that is, \eqref{samantha} holds). The claim on the norms is proved as follows: we have
$$\|A\|=\|T^*T\|\le \|T\|^2$$
on the one hand and
\[
\|A\|\ge\sup_{\|b\|=1}\langle Ab,b\rangle_{{\mathcal B}} =\sup_{\|b\|=1} \langle T^*Tb,b\rangle_{{\mathcal B}}
=\sup_{\|b\|=1}\langle Tb,Tb\rangle_{{\mathcal H}_A} =\|T\|^2,
\]
that is, $\|A\|\ge\|T\|^2$ on the other hand.
Combining the two inequalities we obtain \eqref{bagdad}. We now prove
\eqref{sof-sof}. We have:
\[
\begin{split}
|\langle Ab,c\rangle_{\mathcal B}|&=\langle Tb,Tc\rangle_{\mathcal H}\\
&\le \langle Tb,Tb\rangle_{\mathcal H}^{1/2}
\langle Tc,Tc\rangle_{\mathcal H}^{1/2}\\
&=\langle Ab,b\rangle_{\mathcal B}\langle Ac,c\rangle_{\mathcal B}.
\end{split}
\]
\end{proof}

We will say that $A\leq B$ if $B-A\geq 0.$ Note that
\begin{equation}
\label{etoile}
A\leq B \implies \|A\|\leq\|B\|.
\end{equation}
 Indeed, from \eqref{bagdad} we have:
$$\|A\|=\sup_{\|b\|=1}\langle Ab,b
\rangle_{{\mathcal B}}\le \sup_{\|b\|=1}\langle Bb,b \rangle_{{\mathcal B}}=\|B\|.$$

\begin{definition}\label{defrk}
 Let ${\mathcal H}$ be a Hilbert space of ${\mathcal B}^*$-valued functions
defined on a set $\Omega$ and let $K(z,w)$ be an ${\mathbf L}({\mathcal B},
{\mathcal B}^*)$-valued kernel
defined on $\Omega\times\Omega.$ The kernel $K(z,w)$ is called the reproducing
 kernel of the Hilbert space
${\mathcal H}$ if for every $w\in\Omega$ and $b\in{\mathcal B}$
$K(\cdot,w)b\in\mathcal{H}$ and
\[
\langle f, K(\cdot,w)b\rangle_{\mathcal H}= \langle f(w),
b\rangle_{{\mathcal B}},\quad\forall f\in\mathcal{H}.
\]
\end{definition}

\begin{definition}
Let $K(z,w)$ be an ${\mathbf L}({\mathcal B},{\mathcal B}^*)$-valued
kernel defined on $\Omega\times\Omega.$ The
kernel $K(z,w)$ is said to be positive if for any choice of $z_1,\dots, z_n\in\Omega$ and
$b_1,\dots,b_n\in\mathcal{B}$ it holds that
\[\sum_{j=1}^n\langle K(z_i,z_j)b_j,b_i\rangle_{\mathcal{B}}\geq 0.\]
\end{definition}
\begin{proposition}
\label{p1}
The reproducing kernel $K(z,w)$ of a Hilbert space of ${\mathcal B}^*$-valued functions,
when it exists, is unique and positive.
\end{proposition}

\begin{proposition}
\label{p2}
Let
$K(z,w)$ be an ${\mathbf L}({\mathcal B},{\mathcal B}^*)$-valued
positive kernel defined on
$\Omega\times\Omega.$ Then there exists a unique Hilbert space of
${\mathcal B}^*$-valued functions defined on
$\Omega$ with the reproducing kernel $K(z,w)$.
\end{proposition}

The proofs of these propositions are the same as in the Hilbert space case
and are therefore omitted.

\begin{remark}\label{brem}
One can derive the notion of a reproducing kernel Hilbert space of
 ${\mathcal B}$-valued functions from
Definition \ref{defrk} above, using the natural injection
$\tau$ from ${\mathcal B}$ into
${\mathcal B}^{**}$ defined by \eqref{bastille}.
\end{remark}

\section{Stieltjes integral}
 In this section we define the Stieltjes
integral of a scalar function with respect to a ${\mathbf
L}({\mathcal B},{\mathcal B}^*)$--valued positive measure. We
here follows the analysis presented in \cite[\S 4 p. 19]{MR48:904} for
the case of operators in Hilbert spaces. We consider a separable
Banach space ${\mathcal B}$ and an increasing positive function
$$M\,:\, [a,b]\longrightarrow {\mathbf
L}({\mathcal B},{\mathcal B}^*).$$
Thus, $M(t)\geq 0$ for all $t\in[a,b]$ and moreover
$$
a\le t_1\le t_2\le b\Longrightarrow M(t_2)-M(t_1)\ge 0.$$
Let $f(t)$ be a scalar continuous function on $[a,b]$ and let
$$a=t_0\le \xi_1\le t_1\le \xi_2\le t_2\le\cdots\le \xi_{m}\le
t_m=b$$ be a subdivision of $[a,b]$.
 The Stieltjes integral $\int_a^b f(t)dM(t)$
is defined to be the limit (in the ${\bf L}({\mathcal B},{\mathcal B}^*)$ topology) of the
 sums of the form
$$\sum_{j=1}^m f(\xi_j)(M(t_j)-M(t_{j-1}))$$
as $\sup_{j}|t_j-t_{j-1}|$ goes to $0$.

\begin{theorem}
\label{Sevres-Babylone}
The integral $\int_a^bf(t)dM(t)$ exists.
\end{theorem}
The proof of Theorem \ref{Sevres-Babylone}
is done along the lines of
\cite{MR48:904}. First we need the following lemma.

\begin{lemma}\label{brod}
Let $\alpha_j$ and $\beta_j$ be complex numbers such that
$|\alpha_j|\le |\beta_j|$ ($j=1,2,\ldots m$). Let $H_1,\ldots,
H_m$ be positive operators from ${\mathcal B}$ into ${\mathcal
B}^*$. Then it holds that
\begin{equation*}
\|\sum_{j=1}^m\alpha_j H_j\|\le \|\sum_{j=1}^m |\beta_j|H_j\|.
\end{equation*}
\end{lemma}
\begin{proof} For each $j$ we write $H_j=T_j^*T_j$ where $T_j$ is a
bounded operator from ${\mathcal B}$ into some
Hilbert space ${\mathcal H}_j$ as in Theorem \ref{rosh-hashana}. Then
for $b,c\in{\mathcal B}$ of modulus $1$ we
have:
\[
\begin{split}
|\langle \sum_{j=1}^m\alpha_j H_j b,c\rangle_{{\mathcal B}}| &=
|\sum_{j=1}^m\langle \alpha_j H_j b,c\rangle_{{\mathcal B}}|\\
&\le
\sum_{j=1}^m|\alpha_j| \cdot|\langle T_jb,T_jc\rangle_{{\mathcal H}_j}|\\
&\le \sum_{j=1}^m\sqrt{|\beta_j|} \|T_jb\|_{{\mathcal H}_j}
\sqrt{|\beta_j|}\|T_jc\|_{{\mathcal H}_j}\\
&\le \left(\sum_{j=1}^m|\beta_j|\|T_jb\|^2_{{\mathcal
H}_j}\right)^{1/2}
\left(\sum_{j=1}^m|\beta_j|\|T_jc\|^2_{{\mathcal H}_j}\right)^{1/2}\\
& \le \left(\sum_{j=1}^m|\beta_j|\langle H_jb,
b\rangle_{{\mathcal B}}\right)^{1/2}
\left(\sum_{j=1}^m|\beta_j|\langle H_jc,
c\rangle_{{\mathcal B}}\right)^{1/2}\\
&\le \|\sum_{j=1}^m|\beta_j| H_j\|
\end{split}
\]
where we have used \eqref{bagdad} to get the last inequality.
Thus, taking the supremum on $c$ (of unit norm) we have
$\|\sum_{j=1}^m\alpha_jH_jb\|\le\|\sum_{j=1}^m|\beta_j| H_j\|$,
and hence the required inequality.
\end{proof}

\begin{proof}[Proof of Theorem \ref{Sevres-Babylone}]
It suffices to show that for every $\epsilon>0$ there exists $\delta>0$ such that if
\begin{equation}\label{subdiv}a=t_0\leq t_1\leq\dots\leq t_m=b\end{equation} is a
subdivision of $[a,b]$ such that
$\max_{j}|t_j-t_{j-1}|\leq \delta$ and
$$a=t_0=t_1^0\leq  \dots\leq t_1^{k_1}=t_1=t_2^0\leq\dots\leq t_{m}^{k_{m}} =t_m=b$$
is a continuation of the subdivision \eqref{subdiv} then
for every choice of $\xi_j\in[t_{j-1},t_j]$ and
$\xi_j^\ell\in[t_j^{\ell-1},t_j^\ell]$ it holds that
$$\left\|\sum_{j=1}^m f(\xi_j)(M(t_j)-M(t_{j-1}))-\sum_{j=1}^m
\sum_{\ell=1}^{k_j}
f(\xi_j^\ell)(M(t_j^\ell)-M(t_j^{\ell-1}))\right\|\leq\epsilon.$$
Let us take $\delta$ such that
$$|t^\prime-t^{\prime\prime}|\leq\delta\implies
|f(t^\prime)-f(t^{\prime\prime})|\leq\dfrac{\epsilon}{\|M(b)-M(a)\|}.$$
Then the desired conclusion follows from Lemma \ref{brod}.
\end{proof}
\section{Helly's theorem}
The following theorem is proved in the case of {\sl separable} Hilbert space
in \cite{MR48:904} (see Theorem 4.4 p. 22 there). The proof goes in the same
way in the case of separable Banach spaces. We quote it in a version
adapted to the present setting.
\begin{theorem}
\label{helly}
Let $F_n(t)$ ($t\in[0,2\pi]$) be a sequence of positive increasing
${\mathbf L}({\mathcal B},{\mathcal B}^*)$--valued functions such that
$$F_n(t)\le F_0,\quad n=0,1,\ldots\quad\text{and}\quad t\in[0,2\pi],$$
where $F_0$ is some bounded positive
operator. Then, there exists a subsequence of $F_n$
(which we still denote by $F_n$)
which converges weakly for every $t\in[0,2\pi]$. Moreover, for $f(t)$ a
continuous scalar function we have (in the weak sense, and via the
subsequence):
$$
\int_0^{2\pi}f(t)dF(t)=\lim_{n\rightarrow \infty}\int_0^{2\pi}f(t)dF_n(t).
$$
\end{theorem}
The proof given in \cite{MR48:904} relies on the hypothesis of separability
and on the inequalities
\begin{equation}
\label{rastignac}
\begin{split}
|\langle F_n(t)x,y\rangle_{\mathcal B}|&\le\|F_0\|\cdot\|x\|\cdot\|y\|,\quad
x,y\in{\mathcal B}\\
\sum_{\ell=1}^m |
\langle \Delta_{\ell,n}Fx,y\rangle_{\mathcal B}|&\le 2\|F_0\|\cdot
\|x\|\cdot\|y\|
\end{split}
\end{equation}
where $0=t_0\le t_1<\cdots <t_m=2\pi$ and
$\Delta_{\ell,n}F=F_n(t_\ell)-F_n(t_{\ell-1})$.
The first inequality follows from \eqref{etoile}. We prove the second
one using the factorization given in
Theorem \ref{rosh-hashana}. Using this theorem we write:
$$\Delta_{\ell,n}=T_{\ell,n}^*T_{\ell,n}$$
where $T_{\ell,n}$ is a bounded operator from some Hilbert space
${\mathcal H}_{\ell,n}$ into ${\mathcal B}.$
Then, using \eqref{sof-sof} we have:
\[
\begin{split}
\sum_{\ell =1}^m
|\langle\Delta F_{\ell,n}x,y\rangle_{\mathcal B}|&\le
\sum_{\ell=1}^m
\langle\Delta F_{\ell,n}x,x\rangle_{\mathcal B}^{1/2}
\langle\Delta F_{\ell,n}y,y\rangle_{\mathcal B}^{1/2}\\
&\le
\left(\sum_{\ell=1}^m
\langle\Delta F_{\ell,n}x,x\rangle_{\mathcal B}\right)^{1/2}
\left(\sum_{\ell=1}^m
\langle\Delta F_{\ell,n}y,y\rangle_{\mathcal B}\right)^{1/2}\\
&=\langle (F(2\pi)-F(0))x,x\rangle^{1/2}\langle (F(2\pi)-F(0))y,y
\rangle^{1/2}\\
&\le \langle 2F_0x,x\rangle_{\mathcal B}^{1/2}
\langle 2F_0y,y\rangle_{\mathcal B}^{1/2}\\
&\le 2\|F_0\|\cdot \|x\|\cdot\|y\|.
\end{split}
\]
The proof then proceeds as follows (see
\cite[p. 22]{MR48:904}). One applies the first inequality in \eqref{rastignac}
for $x,y$ in a dense countable set $E$ of ${\mathcal B}$. The functions
$t\mapsto \langle F_n(t)x,y\rangle_{\mathcal B}$ are of bounded variation.
An application of the scalar case of Helly's theorem and
the diagonal process allows to find a subsequence of $F_n$ such that
for all $x,y\in E$ and every $t\in[0,2\pi]$ the limit
$$\lim_{n\rightarrow\infty}\langle F_n(t)x,y\rangle_{\mathcal B}$$
exists. We refer the reader to \cite{MR48:904} for more details.

\section{${\mathbf L}({\mathcal B},{\mathcal B}^*)$-valued
Carath\'eodory functions}
\begin{definition}
Let $\phi(z)$ be an ${\mathbf L}({\mathcal B},{\mathcal B}^*)$-valued function,
  weakly continuous at the origin
in the sense that
\begin{equation}\label{weak} \langle\phi(z)b,b\rangle_{\mathcal{B}}\rightarrow
\langle\phi(0)b,b\rangle_{\mathcal{B}}\text{ as }z\rightarrow 0,\quad\forall b\in\mathcal{B}.
\end{equation}
For a Carath\'eodory function $\phi(z)$ we shall denote by
$\mathcal{L}(\phi)$ the Hilbert space of $\mathcal{B}^*$-valued
functions with the reproducing kernel $k_\phi(z,w)$.
\end{definition}
 We give two representation theorems for Carath\'eodory
functions. In the first we make no assumption on the space ${\mathcal B}$.
Following arguments of Krein
and Langer (see \cite{kl1}), we prove the existence of a realization of the
form \eqref{anaelle2}. The second theorem assumes that
the space ${\mathcal B}$ is separable. We prove that in this case
 the Carath\'eodory functions can be characterized as  functions
analytic in the open unit disk with positive real part. Then we
derive a Herglotz-type representation formula.

\begin{theorem} Let $\Omega$ be a neighborhood of the origin and
let $\phi(z)$ be an ${\mathbf L}({\mathcal B},{\mathcal B}^*)$-valued function
defined in $\Omega$ and weakly continuous at
the origin (in the sense \eqref{weak}).
Then $\phi(z)$ is a Carath\'eodory function if and only if  it admits
the representation
\begin{equation*}
\phi(z)^*\big|_{\mathcal B}=D+C^*(I-z^*V)^{-1}(I+z^*V)C,\quad z\in\Omega,
\end{equation*}
or equivalently,
\begin{equation}
\label{sarah}
\phi(z)=D^*\big|_{\mathcal B}+C^*(I+zV^*)(I-zV^*)^{-1}C,\quad z\in\Omega,
\end{equation}
where $V$ is an isometric operator in some Hilbert space
${\mathcal H}$,  $C$ is a bounded operator from ${\mathcal
B}$ into ${\mathcal H}$ and  $D$ is a purely imaginary
operator from ${\mathcal B}$ into ${\mathcal B}^{*}$ in the sense that
\begin{equation}
\label{D1}
D+D^*\big|_{{\mathcal B}}=0.
\end{equation}
In particular, every Carath\'eodory function has an analytic extension to
the whole open unit disk.
\label{Ivry-sur-Seine}
\end{theorem}

\begin{proof}
Let $\phi(z)$ be a Carath\'eodory function.  First we observe that elements of
$\mathcal{L}(\phi)$ are weakly continuous at the origin:
\[\langle f(w),b\rangle_{\mathcal{B}}\rightarrow
\langle f(0),b\rangle_{\mathcal{B}}\text{ as }w\rightarrow 0,\quad
\forall\,f\in\mathcal{L}(\phi),b\in\mathcal{B}.\]
Indeed, this is a consequence of the Cauchy -- Schwarz inequality as
\[\langle f(w),b\rangle_{\mathcal{B}}-
\langle f(0),b\rangle_{\mathcal{B}}=
\langle f,(k_\phi(\cdot,w)-k_\phi(\cdot,0))b\rangle_{\mathcal{L}(\phi)}\]
and
\[
\|(k_\phi(\cdot, w)-k_\phi(\cdot,0))b\|^2_{\mathcal{L}(\phi)}
=\frac{|w|^2}{1-|w|^2}\Re \langle \phi(w)b,b\rangle_{{\mathcal B}}.
\]

We consider in ${\mathcal L}(\phi)\times{\mathcal L}(\phi)$ the linear
relation $R$ spanned by the pairs
\[
R=
\left(\begin{array}{cc}
\sum k_\phi(z,w_i)w_i^*b_i,\sum k_\phi(z,w_i)b_i-k_\phi(z,0)(\sum b_i)
\end{array}\right)
\]
where all the $b_i\in{\mathcal B}$, the $w_i\in\Omega$ and where all the
sums are finite. This
relation is densely defined because of the weak continuity of the elements of
${\mathcal L}(\phi)$ at the origin.
Indeed, let $f\in{\mathcal L}(\phi)$ be orthogonal to the domain of $R$. Then,
\[
\langle f, k_\phi(\cdot,w)b\rangle_{{\mathcal L}(\phi)}=\langle f(w),b
\rangle_{{\mathcal B}}=0\]
for all $b\in{\mathcal B}$ and all
points $w\not =0$ in the domain of $f$. Thus $f(w)=0$ at all these points $w$
and the continuity hypothesis implies that also $f(0)=0$.
The relation $R$ is readily seen to be isometric. Its closure is thus the
graph of an isometry, which we call $V$. We have:
\[
V(k_\phi(z,w)w^*b)=k_\phi(z,w)b-k_\phi(z,0)b,\]
and in particular
\begin{equation}
\label{Alfort-Ecole-veterinaire}
(I-w^*V)^{-1}k_\phi(\cdot, 0)b=k_\phi(\cdot, w)b.
\end{equation}
Denote by $C$ the map
$$C\,:\,{\mathcal B}\longrightarrow {\mathcal L}(\phi),\quad(Cb)(z)
\defi k_\phi(z,0)b.$$
Then, for $f\in{\mathcal L}(\phi)$, $C^*f=f(0)$. 
Applying $C$ on the left on both sides of
\eqref{Alfort-Ecole-veterinaire} we obtain
\[\frac{\phi(0)+\phi(w)^*\big|_{\mathcal B}}{2}b=C^*(I-w^*V)^{-1}Cb.\]
Since $$C^*Cb=\frac{\phi(0)+\phi(0)^*\big|_{\mathcal B}}{2}b$$
we obtain
\[
\begin{split}
\phi(0)+\phi(w)^*\big|_{\mathcal B}&=2C^*(I-w^*V)^{-1}C-C^*C+C^*C\\
                                   &=C^*(I-w^*V)^{-1}(I+w^*V)C+C^*C
\end{split}
\]
so that
\[
\phi(w)^*\big|_{{\mathcal B}}+\frac{\phi(0)-\phi(0)^*\big|_{\mathcal B}}{2}
=
C^*(I-w^*V)^{-1}(I+w^*V)C,
\]
which gives the required formula with
\begin{equation}
\label{D} D=\frac{\phi(0)-\phi(0)^*\big|_{\mathcal B}}{2}.
\end{equation}

We now prove the converse statement and first compute
$$\langle k_\phi(z,w)x,y\rangle_{{\mathcal B}}
\quad{\rm for}\quad x,y\in{\mathcal B}.$$
We have
\[
\begin{split}
\langle\phi(w)^*\big|_{\mathcal B}x,y\rangle_{{\mathcal B}}&=
\langle Dx,y\rangle_{{\mathcal B}}+\langle
(I-w^*V)^{-1}(I+w^*V)Cx,Cy\rangle_{{\mathcal L}(\phi)}.
\end{split}
\]

We have, with $\tau$ the natural injection from ${\mathcal B}$ into
${\mathcal B}^{**}$ (see \eqref{bastille}):
\begin{equation}
\begin{split}
\langle \phi(z)x,y\rangle_{{\mathcal B}} &= \overline{\langle \tau
y ,\phi(z)x\rangle_{{\mathcal B}^*}}\\
&=\overline{\langle \phi(z)^*\tau y,x\rangle_{{\mathcal B}^*}}\\
&= \overline{\langle \phi(z)^*\big|_{\mathcal B}y,x
\rangle_{{\mathcal B}}}\\
&=\overline{\langle D y,x\rangle_{{\mathcal B}}}+
\langle(I-zV^*)^{-1}(I+zV^*)Cx,Cy\rangle_{{\mathcal L}(\phi)}\\
&=\langle D^*\big|_{\mathcal B}x ,y\rangle_{{\mathcal B}}+
\langle(I-zV^*)^{-1}(I+zV^*)Cx,Cy\rangle_{{\mathcal L}(\phi)},
\end{split}
\label{elinor}
\end{equation}
and so
$$\langle k_\phi(z,w)x,y\rangle_{{\mathcal B}}
=\langle (I-zV^*)^{-1}Cx,(I-wV^*)^{-1}Cy\rangle_{{\mathcal L}(\phi)},$$
from which follows the positivity of $k_\phi(z,w)$. Finally,
\eqref{elinor} also implies \eqref{sarah} and this concludes the proof.
\end{proof}

\begin{remark}
Although the above argument is very close to the one in
\cite[p. 365--366]{kl1} we note the following: we use a concrete
space (the space ${\mathcal L}(\phi)$) to build the relation rather than
abstract elements and the relation $R$ is defined slightly differently.
\end{remark}
\begin{remark}\label{indef}
As already mentioned, the above argument still goes through when the
kernel has a finite number of negative squares. In this case the space
$\mathcal{L}(\phi)$ is a Pontryagin space.  For $b\in
{\mathcal B}$ and sufficiently small $h\in{\mathbb C}$ we consider
the functions $f_h(z)=(K(z,w+h)-K(z,w))b$, which  have the following
properties:
\begin{align}
\nonumber
\lim_{h\rightarrow 0}\langle f_h,f_h\rangle_{{\mathcal L}(\phi)}&=0,\\
\lim_{h\rightarrow 0}\langle f,f_h\rangle_{{\mathcal L}(\phi)}&=0,
\quad\forall f\in\spa\{k_\phi(\cdot,w)b\}.
\nonumber
\end{align}
It follows from the convergence criteria in Pontryagin spaces (see
\cite{ikl}, \cite[p. 356]{kl1}) that
$$\lim_{h\rightarrow 0}\langle f,f_h\rangle_{{\mathcal L}(\phi)}=0,\quad
\forall f\in\mathcal{L}(\phi).$$ The fact that the relation is the graph of
an isometric operator is
proved in \cite[Theorem 1.4.2 p. 29]{adrs}. This follows from a theorem of
Shmulyan which states that a
contractive relation between Pontryagin spaces of same index is the graph
 of a contractive operator; see
\cite{s} and \cite[Theorem 1.4.1 p. 27]{adrs}.
\end{remark}

\begin{theorem}
Let ${\mathcal B}$ be a separable Banach space and let $\phi(z)$ be a
 ${\mathbf L}({\mathcal B},{\mathcal
B}^*)$--valued function analytic in the open unit disk, such that
\[\phi(z)+{\phi(z)^*}\bigm|_{ \mathcal{B}}\,
\geq 0.\] Then there exists an increasing ${\mathbf L}({\mathcal B},{\mathcal B}^*)$--valued
function $M(t)$
($t\in[0,2\pi]$) and a purely imaginary operator $D$ (that is, satisfying  \eqref{D1})
such that
\[
\phi(z)=D+\int_0^{2\pi}\frac{e^{it}+z}{e^{it}-z}dM(t),
\]
where the integral is defined in the weak sense. Furthermore the
kernel $k_\phi(z,w)$ is positive in ${\mathbb D}$.
\label{Antony}
\end{theorem}
\begin{proof} We follow the arguments in \cite{MR48:904}, and will apply
Theorem \ref{helly}. The separability hypothesis of ${\mathcal B}$ is used at
this point.\\

We first assume that $\phi(z)$
is analytic in $|z|<1+\epsilon$ with $\epsilon>0$. We have
(the existence of the integrals follows from Theorem
\ref{Sevres-Babylone}):
\[
\begin{split}
\frac{1}{4\pi}
\int_0^{2\pi}\phi(e^{it})\frac{e^{it}+z}{e^{it}-z}dt&=\phi(z)-\frac{\phi(0)}{2}\\
\frac{1}{4\pi}\int_0^{2\pi}(\phi(e^{it})^*\big|_{\mathcal B}
\frac{e^{it}+z}{e^{it}-z}dt&= \frac{1}{2\pi
}\int_0^{2\pi}(\phi(e^{-it})^*\big|_{\mathcal B}
\frac{e^{-it}+z}{e^{-it}-z}\\
&=\frac{(\phi(0))^*\big|_{\mathcal B}}{2}.
\end{split}
\]
Thus,
\[
\phi(z)=D+\int_0^{2\pi}\frac{\left(\phi(e^{it})+(\phi(e^{it}))^*\big|_{\mathcal
B}\right)}{4\pi}\frac{e^{it}+z}{e^{it}-z}dt,
\] with $D$ as in \eqref{D} and the formula for general
$\phi(z)$ follows from Helly's theorem applied to the measures
$$
\frac{\left(\phi(re^{it})+(\phi(re^{it}))^*\big|_{\mathcal
B}\right)}{4\pi}dt,\quad r<1$$
(or more precisely to a sequence $r_n\rightarrow 1$).\\

We now prove the positivity of the kernel $k_\phi(z,w)$ and first
assume that $\phi(z)$ is analytic in $|z|<1+\epsilon$ as above.
We have:
\[
\frac{1}{4\pi}\int_0^{2\pi}(\phi(e^{-it}))^*\big|_{\mathcal B}
\frac{e^{it}+z}{e^{it}-z}dt=(\phi(z^*))^*\big|_{\mathcal B}
-\frac{(\phi(0))^*\big|_{\mathcal B}}{2}.
\]
Thus
\[
(\phi(z^*))^*\big|_{\mathcal B} -\frac{(\phi(0))^*\big|_{\mathcal
B}}{2}= \frac{1}{4\pi}\int_0^{2\pi}(\phi(e^{it}))^*\big|_{\mathcal
B} \frac{e^{-it}+z}{e^{-it}-z}dt
\]
and so
$$
(\phi(z^*))^*\big|_{\mathcal B}=
-D+\int_0^{2\pi}\frac{\left(\phi(e^{it})+(\phi(e^{it}))^*\big|_{\mathcal
B}\right)}{4\pi}\frac{1+ze^{it}}{1-ze^{it}}dt
$$
since
$$
\frac{1}{4\pi}\int_0^{2\pi}\phi(e^{it})\frac{1+ze^{it}}{1-ze^{it}}dt
=
\frac{1}{4\pi}\int_0^{2\pi}\phi(e^{-it})\frac{e^{it}+z}{e^{it}-z}dt=\frac{\phi(0)}{2}.$$
Thus,
$$k_\phi(z,w)=\int_0^{2\pi}
\frac{\left(\phi(e^{it})+(\phi(e^{it}))^*\big|_{\mathcal B}\right)}{4\pi}\frac{1}{(e^{it}-z)(e^{it}-w)^*}dt.$$
The positivity follows. The case of general $\phi(z)$ is done by approximation using Helly's theorem; indeed
using Theorem \ref{helly} we have for general $\phi(z)$:
\[
\begin{split}
\langle k_\phi(w_\ell,w_j)b_j,b_\ell\rangle_{\mathcal B}&=
\\
&\hspace{-2cm}=
\lim_{r\rightarrow 1}\int_0^{2\pi}
\left\langle
\frac{\left(\phi(re^{it})+(\phi(re^{it}))^*\big|_{\mathcal
B}\right)}{4\pi}\frac{1}{(e^{it}-w_\ell)(e^{it}-w_j)^*}b_j,b_\ell
\right\rangle_{\mathcal B}dt\\
&\hspace{-2cm}=
\lim_{r\rightarrow 1}\int_0^{2\pi}
\left\langle
\frac{\left(\phi(re^{it})+(\phi(re^{it}))^*\big|_{\mathcal
B}\right)}{4\pi}\frac{b_j}{(e^{it}-w_j)^*},\frac{b_\ell}{(e^{it}-w_\ell)^*}
\right\rangle_{\mathcal B}dt.
\end{split}
\]
\end{proof}

\section{The case of ${\mathbf L}({\mathcal B}^{*},{\mathcal B})$--valued
functions} We turn to the case of ${\mathbf L}({\mathcal
B}^{*},{\mathcal B})$--valued  functions. 
Using the natural
injection $\tau$
$$\mathcal{B}\stackrel{\tau}{\mapsto}\mathcal{B}^{**}$$
defined by \eqref{bastille} we shall say that an ${\mathbf
L}({\mathcal B}^{*},{\mathcal B})$--valued function $\phi(z)$ is
a Carath\'eodory function if the ${\mathbf L}({\mathcal
B}^{*},{\mathcal B}^{**})$-valued function $\tau\phi(z)$ is a
Carath\'eodory function.

\begin{theorem} An  ${\mathbf L}({\mathcal B}^*,{\mathcal B})$--valued
function $\phi(z)$ defined in a neighborhood of the origin
and weakly continuous at the origin is a Carath\'eodory function
if and only if it admits the representation
\begin{equation*}
\phi(z)^*=D+C^*(I-z^*V)^{-1}(I+z^*V)C,
\end{equation*}
or, equivalently,
\[
\tau\phi(z)=D^*\big|_{{\mathcal B}_*}+C^*(I+zV^*)(I-zV^*)^{-1}C,\]
where $V$ is an isometric operator in some Hilbert space
${\mathcal H}$, where $C$ is a bounded operator from ${\mathcal
B}^*$ into ${\mathcal H}$ and where $D$ is a purely imaginary
operator from ${\mathcal B}^*$ into ${\mathcal B}^{**}$.

In particular $\phi(z)$ has an analytic extension to the whole open unit disk.
\end{theorem}

\begin{proof} By Theorem \ref{Ivry-sur-Seine} (with $\mathcal{B}$ replaced by
$\mathcal{B}^*$),
$\phi(z)$ is a Carath\'eodory function if and only if
\[(\tau\phi(z))^*\bigm|_{\mathcal{B}^*}=D+C^*(I-z^*V)^{-1}(I+z^*V)C,\]
where $C, V, D$ have the stated properties.
But $(\tau\phi(z))^*\bigm|_{\mathcal{B}^*}=\phi(z)^*.$
\end{proof}

\bibliographystyle{plain}

\begin{thebibliography}{10}

\bibitem{MR2003f:47021}
D.~Alpay, V.~Bolotnikov, A.~Dijksma, and B.~Freydin.
\newblock Nonstationary analogues of the {H}erglotz representation theorem for
  unbounded operators.
\newblock {\em Arch. Math. (Basel)}, 78(6):465--474, 2002.

\bibitem{MR1704663}
D.~Alpay, A.~Dijksma, and Y.~Peretz.
\newblock Nonstationary analogs of the {H}erglotz representation theorem: the
  discrete case.
\newblock {\em J. Funct. Anal.}, 166(1):85--129, 1999.

\bibitem{adrs}
D.~Alpay, A.~Dijksma, J.~Rovnyak, and H.~de~Snoo.
\newblock {\em {Schur} functions, operator colligations, and reproducing kernel
  {P}ontryagin spaces}, volume~96 of {\em Operator theory: {A}dvances and
  {A}pplications}.
\newblock Birkh{\" a}user Verlag, Basel, 1997.

\bibitem{ad1}
D.~Alpay and H.~Dym.
\newblock Hilbert spaces of analytic functions, inverse scattering and operator
  models, {I}.
\newblock {\em Integral Equation and Operator Theory}, 7:589--641, 1984.

\bibitem{ad2}
D.~Alpay and H.~Dym.
\newblock Hilbert spaces of analytic functions, inverse scattering and operator
  models, {II}.
\newblock {\em Integral Equation and Operator Theory}, 8:145--180, 1985.

\bibitem{MR1960423}
D.~Alpay and I.~Gohberg.
\newblock A trace formula for canonical differential expressions.
\newblock {\em J. Funct. Anal.}, 197(2):489--525, 2003.


\bibitem{dbr1}
{L. de} Branges and J.~Rovnyak.
\newblock Canonical models in quantum scattering theory.
\newblock In C.~Wilcox, editor, {\em Perturbation theory and its applications
  in quantum mechanics}, pages 295--392. Wiley, {N}ew {Y}ork, 1966.

\bibitem{dbr2}
{L. de} Branges and J.~Rovnyak.
\newblock {\em Square summable power series}.
\newblock {Holt, Rinehart and Winston, New {Y}ork}, 1966.

\bibitem{dbs}
{L. de} Branges and L.A. Shulman.
\newblock Perturbation theory of unitary operators.
\newblock {\em J. Math. Anal. Appl.}, 23:294--326, 1968.

\bibitem{MR48:904}
M.~S. Brodski{\u\i}.
\newblock {\em Triangular and {J}ordan representations of linear operators}.
\newblock American Mathematical Society, Providence, R.I., 1971.
\newblock Translated from the Russian by J. M. Danskin, Translations of
  Mathematical Monographs, Vol. 32.

\bibitem{Dmk}
H.~Dym and H.P. McKean.
\newblock {\em Gaussian processes, function theory and the inverse spectral
  problem}.
\newblock Academic {P}ress, 1976.

\bibitem{MR1973084}
B.~Farkas and M.~Matolcsi.
\newblock Positive forms on {B}anach spaces.
\newblock {\em Acta Math. Hungar.}, 99(1-2):43--55, 2003.

\bibitem{MR1821917}
F.~Gesztesy, N.~Kalton, K.A. Makarov, and E.~Tsekanovskii.
\newblock Some applications of operator-valued {H}erglotz functions.
\newblock In D.~Alpay and V.~Vinnikov, editors, {\em Operator theory, system
  theory and related topics (Beer-Sheva/Rehovot, 1997)}, volume 123 of {\em
  Oper. Theory Adv. Appl.}, pages 271--321. Birkh\"auser, Basel, 2001.

\bibitem{MR1759548}
F.~Gesztesy and K.A. Makarov.
\newblock Some applications of the spectral shift operator.
\newblock In {\em Operator theory and its applications (Winnipeg, MB, 1998)},
  pages 267--292. Amer. Math. Soc., Providence, RI, 2000.

\bibitem{MR647140}
J.~G{\'o}rniak and A.~Weron.
\newblock Aronszajn-{K}olmogorov type theorems for positive definite kernels in
  locally convex spaces.
\newblock {\em Studia Math.}, 69(3):235--246, 1980/81.

\bibitem{ikl}
I.S. Iohvidov, M.G. Kre\u{\i}n, and H.~Langer.
\newblock {\em Introduction to the spectral theory of operators in spaces with
  an indefinite metric}.
\newblock Akademie--Verlag, Berlin, 1982.

\bibitem{kl1}
M.G. Kre{\u\i}n and H.~Langer.
\newblock {\"{U}}ber die verallgemeinerten {R}esolventen und die
  charakteristische {F}unktion eines isometrischen {O}perators im {R}aume ${\Pi
  _k}$.
\newblock In {\em Hilbert space operators and operator algebras (Proc. Int.
  Conf. Tihany, 1970)}, pages 353--399. North--Holland, Amsterdam, 1972.
\newblock Colloquia Math. Soc. J\'{a}nos Bolyai.

\bibitem{masani78}
P.~Masani.
\newblock Dilations as propagators of hilbertian varieties.
\newblock {\em SIAM J. Math. Anal.}, 9:414--456, 1978.

\bibitem{MR899274}
S.~N. Naboko.
\newblock On the boundary values of analytic operator-valued functions with a
  positive imaginary part.
\newblock {\em Zap. Nauchn. Sem. Leningrad. Otdel. Mat. Inst. Steklov. (LOMI)},
  157(Issled. po Linein. Operator. i Teorii Funktsii. XVI):55--69, 179, 1987.

\bibitem{MR1036844}
S.~N. Naboko.
\newblock Nontangential boundary values of operator {$R$}-functions in a
  half-plane.
\newblock {\em Algebra i Analiz}, 1(5):197--222, 1989.

\bibitem{pedrick}
G.~Pedrick.
\newblock Theory of reproducing kernels for {H}ilbert spaces of vector valued
  functions.
\newblock Studies in eigenvalues problems~19, University of Kansas, Department
  of {M}athematics, July 1957.

\bibitem{MR58:12429a}
M.~Reed and B.~Simon.
\newblock {\em Methods of modern mathematical physics. {I}. {F}unctional
  analysis}.
\newblock Academic Press, New York, 1972.

\bibitem{s}
Y.L. Shmul'yan.
\newblock Division in the class of ${J}$--expansive operators.
\newblock {\em Math. {S}b.}, 116:516---525, 1967.

\bibitem{MR626346}
N.~N. Vakhania.
\newblock {\em Probability distributions on linear spaces}.
\newblock North-Holland Publishing Co., New York, 1981.
\newblock Translated from the Russian by I. I. Kotlarski, North-Holland Series
  in Probability and Applied Mathematics.

\end{thebibliography}

\end{document}